\documentclass[leqno,12pt]{amsart}
\usepackage{amsfonts}
\usepackage{amsmath,amssymb,amsthm}

\newcommand{\R}{\mathbb R}

\renewcommand{\span}{\mathrm{span}}
\newcommand{\tr}{\mathrm{tr}}

\newtheorem{thm}{Theorem}[section]

\theoremstyle{definition}

\theoremstyle{remark}

\newcommand{\ds}{\displaystyle}

\begin{document}

\title[General Rotational Surfaces in the Four-dimensional Minkowski
Space] {General Rotational Surfaces in the Four-dimensional
Minkowski Space}

\author{Georgi Ganchev and Velichka Milousheva}
\address{Institute of Mathematics and Informatics, Bulgarian Academy of Sciences,
Acad. G. Bonchev Str. bl. 8, 1113 Sofia, Bulgaria}
\email{ganchev@math.bas.bg}
\address{Institute of Mathematics and Informatics, Bulgarian Academy of Sciences,
Acad. G. Bonchev Str. bl. 8, 1113, Sofia, Bulgaria;   "L. Karavelov"
Civil Engineering Higher School, 175 Suhodolska Str., 1373 Sofia,
Bulgaria} \email{vmil@math.bas.bg}

\subjclass[2000]{Primary 53A35, Secondary 53B25}
\keywords{Surfaces in the 4-dimensional Minkowski space, general
rotational surfaces, minimal surfaces, flat surfaces, surfaces
with flat normal connection}

\begin{abstract}
General rotational surfaces as a source of examples of surfaces in the four-dimensional Euclidean space have been introduced by C. Moore.
In this paper we consider the analogue of these surfaces in the Minkowski 4-space.
On the base of our invariant theory of spacelike surfaces we study general rotational surfaces with special invariants. We describe
 analytically the flat general rotational surfaces
and the general rotational surfaces with flat normal connection. We classify  completely the minimal general rotational surfaces and
the general rotational surfaces consisting of parabolic points.
\end{abstract}

\maketitle

\section{Introduction}

The local theory of spacelike surfaces in the four-dimensional Minkowski space
$\R^4_1$ was developed by the present authors in \cite {GM5}. Our approach to this theory is based on the
introduction of an invariant linear map of Weingarten-type in the tangent plane at
any point of the surface.
This invariant map allowed  us to introduce principal lines
and a geometrically determined  moving frame field at each point of the surface. Writing derivative formulas
of Frenet-type for this frame field, we obtained eight invariant
functions   $\gamma_1, \, \gamma_2, \, \nu_1,\, \nu_2, \, \lambda, \, \mu,
\, \beta_1, \beta_2$,   and proved a fundamental theorem of Bonnet-type, stating
that these eight invariants under some natural conditions
determine the surface up to a motion in $\R^4_1$.

The basic geometric classes of surfaces in $\R^4_1$  are characterized by conditions on these invariant functions.
For example, surfaces with flat normal connection are characterized by the condition $\nu_1 = \nu_2$,
minimal surfaces are  described by $\nu_1 + \nu_2 = 0$, Chen surfaces are characterized by  $\lambda = 0$.

Rotational surfaces are basic source of examples of many geometric classes of surfaces. In \cite{M} C. Moore introduced general rotational surfaces
in the four-dimensional Euclidean space $\R^4$ and described a special case of general rotational surfaces with constant Gauss curvature \cite{M2}.

In the present paper we consider spacelike  general rotational surfaces which are analogous to the general rotational surfaces of C. Moore.
We apply the invariant theory  of spacelike surfaces in $\R^4_1$ to the class  of general rotational surfaces with plane meridians.
Using the invariants of these surfaces
we describe analytically the flat general rotational surfaces (Theorem \ref{T:Th-flat})
and the general rotational surfaces with flat normal connection (Theorem  \ref{T:Th-flat normal}).
In  Theorem \ref{T:Th-parabolic points} we give the complete classification of general rotational surfaces consisting of parabolic points. The classification of
minimal general rotational surfaces is given in Theorem \ref{T:Th-minimal}.

\section{Preliminaries}

Let $\R^4_1$  be the four-dimensional Minkowski space  endowed with the metric
$\langle , \rangle$ of signature $(3,1)$ and  $Oe_1e_2e_3e_4$ be a
fixed orthonormal coordinate system, i.e. $e_1^2 =
e_2^2 = e_3^2 = 1, \, e_4^2 = -1$, giving the orientation of
$\R^4_1$.
A surface $M^2: z = z(u,v), \, \, (u,v) \in {\mathcal D}$
(${\mathcal D} \subset \R^2$) in $\R^4_1$ is said to be
\emph{spacelike} if $\langle , \rangle$ induces  a Riemannian
metric $g$ on $M^2$. Thus at each point $p$ of a spacelike surface
$M^2$ we have the following decomposition:
$$\R^4_1 = T_pM^2 \oplus N_pM^2$$
with the property that the restriction of the metric $\langle ,
\rangle$ onto the tangent space $T_pM^2$ is of signature $(2,0)$,
and the restriction of the metric $\langle , \rangle$ onto the
normal space $N_pM^2$ is of signature $(1,1)$.

Denote by $\nabla'$ and $\nabla$ the Levi Civita connections on $\R^4_1$ and $M^2$, respectively.
Let $x$ and $y$ be vector fields tangent to $M$ and $\xi$ be a normal vector field.
The formulas of Gauss and Weingarten give a decomposition of the vector fields $\nabla'_xy$ and
$\nabla'_x \xi$ into a tangent and a normal component:
$$\begin{array}{l}
\vspace{2mm}
\nabla'_xy = \nabla_xy + \sigma(x,y);\\
\vspace{2mm}
\nabla'_x \xi = - A_{\xi} x + D_x \xi,
\end{array}$$
which define the second fundamental tensor $\sigma$, the normal connection $D$
and the shape operator $A_{\xi}$ with respect to $\xi$.

The mean curvature vector  field $H$ of $M^2$ is defined as
$H = \ds{\frac{1}{2}\,  \tr\, \sigma}$,
i.e. given a local orthonormal frame $\{x,y\}$ of the tangent bundle,
$H = \ds{\frac{1}{2} \left(\sigma(x,x) +  \sigma(y,y)\right)}$.

Let
$M^2: z=z(u,v), \,\, (u,v) \in \mathcal{D}$ $(\mathcal{D} \subset \R^2)$
be a local parametrization on a
spacelike surface in $\R^4_1$.
The tangent space at an arbitrary point $p=z(u,v)$ of $M^2$ is
$T_pM^2 = \span \{z_u,z_v\}$, where $\langle z_u,z_u \rangle > 0$, $\langle z_v,z_v \rangle > 0$.
We use the standard denotations
$E(u,v)=\langle z_u,z_u \rangle, \; F(u,v)=\langle z_u,z_v
\rangle, \; G(u,v)=\langle z_v,z_v \rangle$ for the coefficients
of the first fundamental form
$$I(\lambda,\mu)= E \lambda^2 + 2F \lambda \mu + G \mu^2,\quad
\lambda, \mu \in \R.$$
 Since $I(\lambda, \mu)$ is positive
definite we set $W=\sqrt{EG-F^2}$.
We choose a normal frame field $\{n_1, n_2\}$ such that $\langle
n_1, n_1 \rangle =1$, $\langle n_2, n_2 \rangle = -1$, and the
quadruple $\{z_u,z_v, n_1, n_2\}$ is positively oriented in
$\R^4_1$.
Then we have the following derivative formulas:
$$\begin{array}{l}
\vspace{2mm} \nabla'_{z_u}z_u=z_{uu} = \Gamma_{11}^1 \, z_u +
\Gamma_{11}^2 \, z_v + c_{11}^1\, n_1 - c_{11}^2\, n_2;\\
\vspace{2mm} \nabla'_{z_u}z_v=z_{uv} = \Gamma_{12}^1 \, z_u +
\Gamma_{12}^2 \, z_v + c_{12}^1\, n_1 - c_{12}^2\, n_2;\\
\vspace{2mm} \nabla'_{z_v}z_v=z_{vv} = \Gamma_{22}^1 \, z_u +
\Gamma_{22}^2 \, z_v + c_{22}^1\, n_1 - c_{22}^2\, n_2,\\
\end{array}$$
where $\Gamma_{ij}^k$ are the Christoffel's symbols and the functions
$c_{ij}^k, \,\, i,j,k = 1,2$ are given by
$$\begin{array}{lll}
\vspace{2mm}
c_{11}^1 = \langle z_{uu}, n_1 \rangle; & \qquad   c_{12}^1 = \langle z_{uv},
n_1 \rangle; & \qquad  c_{22}^1 = \langle z_{vv}, n_1 \rangle;\\
\vspace{2mm}
c_{11}^2 = \langle z_{uu}, n_2 \rangle; & \qquad  c_{12}^2 = \langle z_{uv},
n_2 \rangle; & \qquad c_{22}^2 = \langle z_{vv}, n_2 \rangle.
\end{array} $$

Obviously, the surface $M^2$ lies in a 2-plane if and only if
$M^2$ is totally geodesic, i.e. $c_{ij}^k=0, \; i,j,k = 1, 2.$ So,
we assume that at least one of the coefficients $c_{ij}^k$ is not
zero.

The second fundamental form $II$ of the surface $M^2$ at a point
$p \in M^2$ is introduced by the following functions
\begin{equation} \notag
L = \ds{\frac{2}{W}} \left|%
\begin{array}{cc}
\vspace{2mm}
  c_{11}^1 & c_{12}^1 \\
  c_{11}^2 & c_{12}^2 \\
\end{array}%
\right|; \quad
M = \ds{\frac{1}{W}} \left|%
\begin{array}{cc}
\vspace{2mm}
  c_{11}^1 & c_{22}^1 \\
  c_{11}^2 & c_{22}^2 \\
\end{array}%
\right|; \quad
N = \ds{\frac{2}{W}} \left|%
\begin{array}{cc}
\vspace{2mm}
  c_{12}^1 & c_{22}^1 \\
  c_{12}^2 & c_{22}^2 \\
\end{array}%
\right|.
\end{equation}
Let
$X=\lambda z_u+\mu z_v, \,\, (\lambda,\mu)\neq(0,0)$ be a tangent
vector at a point $p \in M^2$. Then
$$II(\lambda,\mu)=L\lambda^2+2M\lambda\mu+N\mu^2, \quad \lambda, \mu \in {\R}.$$
The second fundamental
form $II$ is invariant up to the orientation of the tangent space
or the normal space of the surface.

The condition $L = M = N = 0$  characterizes points at which
the space  $\{\sigma(x,y):  x, y \in T_pM^2\}$ is one-dimensional.
We call such points  \emph{flat points} of the surface \cite{GM5}.
These points are analogous to flat points in the theory of surfaces in $\R^3$ and $\R^4$ \cite{GM1, GM2}.
In \cite{GM5} we gave a local geometric description of spacelike surfaces consisting of flat points
 proving that
any spacelike surface consisting of flat points whose mean
curvature vector at any point is a non-zero spacelike vector or
timelike vector  either lies in a hyperplane of $\R^4_1$ or is
part of a developable ruled surface in $\R^4_1$.

Further we consider surfaces free of flat points, i.e.  $(L, M, N) \neq (0,0,0)$.

The second fundamental form $II$ determines conjugate, asymptotic, and
principal tangents at a point $p$ of $M^2$ in the standard way.
A line $c: u=u(q), \; v=v(q); \; q\in J \subset \R$ on $M^2$ is
said to be an \emph{asymptotic line}, respectively a
\textit{principal line}, if its tangent at any point is
asymptotic, respectively  principal. The surface $M^2$ is
parameterized by principal lines if and only if $F=0, \,\, M=0.$

The second fundamental form $II$ generates  two invariant functions:
$$k =  \frac{LN - M^2}{EG - F^2}, \qquad
\varkappa = \frac{EN+GL-2FM}{2(EG-F^2)}.$$

The functions $k$ and $\varkappa$
are invariant under changes of the parameters of the surface and
changes of the normal frame field.
 The sign of $k$ is invariant under congruences
 and the sign of $\varkappa$ is invariant under motions
in $\R^4_1$. However, the sign of $\varkappa$ changes under
symmetries with respect to a hyperplane in $\R^4_1$. It turns out that
the invariant $\varkappa$ is the curvature of the normal
connection of the surface
(see \cite{GM5}). The number of asymptotic tangents at a
point of $M^2$ is determined by the sign of the invariant $k$.

As in the theory of surfaces in $\R^3$ and $\R^4$  the invariant
$k$ divides the points of $M^2$ into the following  types: \emph{elliptic} ($k > 0$),
\emph{parabolic} ($k = 0$), and \emph{hyperbolic} ($k < 0$).

Let $H$
be the normal mean curvature vector field. Recall that a
surface $M^2$ is said to be \textit{minimal} if its mean curvature vector vanishes identically, i.e. $H =0$.
 The minimal surfaces are characterized in terms of
the invariants $k$ and $\varkappa$ by the  equality \cite{GM5}
$$\varkappa^2 - k =0.$$

It is interesting to note that the
''\emph{umbilical}'' points, i.e. points at which the coefficients of the
first and the second fundamental forms are proportional
($L = \rho E, \, M = \rho F, \,N = \rho G, \, \rho \neq 0$), are
exactly the points at which the mean curvature vector $H$ is zero.
So, the spacelike surfaces consisting of ''umbilical'' points in
$\R^4_1$ are exactly the minimal surfaces. If $M^2$ is a spacelike surface free of ''umbilical''
 points ($H \neq 0$ at each point), then there exist exactly two principal tangents.

\vskip 2mm
Considering spacelike surfaces in $\R^4_1$ whose mean
curvature vector at any point is a non-zero spacelike vector or
timelike vector, on the base of the principal lines we introduced
a geometrically determined orthonormal frame field $\{x,y,b,l\}$ at each point of such a surface  \cite{GM5}.
The tangent vector fields $x$ and $y$ are collinear with  the principal directions,
the normal vector field $b$ is collinear with the mean curvature vector field $H$.
Writing derivative formulas
of Frenet-type for this frame field, we obtained eight invariant
functions $\gamma_1, \, \gamma_2, \, \nu_1,\, \nu_2, \, \lambda, \, \mu,
\, \beta_1, \beta_2$, which determine the surface up to a rigid motion in $\R^4_1$.

The invariants $\gamma_1, \, \gamma_2, \, \nu_1,\, \nu_2, \, \lambda, \, \mu,
\, \beta_1$, and $\beta_2$ are determined by the geometric frame field $\{x,y,b,l\}$ as follows
$$\nu_1 = \langle \nabla'_xx, b\rangle, \qquad \nu_2 = \langle \nabla'_yy, b\rangle, \qquad \lambda = \langle \nabla'_xy, b\rangle,
\qquad \mu = \langle \nabla'_xy, l\rangle,$$
\begin{equation} \notag
\gamma_1 = - y(\ln \sqrt{E}), \quad  \gamma_2 = - x(\ln \sqrt{G}), \quad \beta_1 = \langle \nabla'_xb, l\rangle, \quad
\beta_2 = \langle \nabla'_yb, l\rangle.
\end{equation}

The invariants $k$, $\varkappa$, and the Gauss curvature $K$ of
$M^2$ are expressed by the functions $\nu_1, \nu_2, \lambda, \mu$
as follows:
\begin{equation} \notag
k = - 4\nu_1\,\nu_2\,\mu^2, \quad \quad \varkappa = (\nu_1-\nu_2)\mu, \quad \quad
K = \varepsilon (\nu_1\,\nu_2- \lambda^2 + \mu^2),
\end{equation}
where $\varepsilon = sign \langle H, H \rangle$.
The  norm $\Vert H \Vert$ of the mean curvature
vector is expressed as
\begin{equation} \notag
\Vert H \Vert = \displaystyle{ \frac{|\nu_1 + \nu_2|}{2} = \frac{\sqrt{\varkappa^2-k}}{2 |\mu |}}.
\end{equation}
If $M^2$ is a spacelike surface whose mean curvature vector at any point is a non-zero spacelike vector or timelike vector, then
$M^2$ is minimal if and only if $\nu_1 + \nu_2 = 0$.

The geometric meaning of the invariant $\lambda$ is connected with the notion of Chen submanifolds.
Let $M$ be an $n$-dimensional submanifold of
$(n+m)$-dimensional Riemannian manifold $\widetilde{M}$ and $\xi$
be a normal vector field of $M$. B.-Y. Chen \cite{Chen1} defined
 the \emph{allied vector field} $a(\xi)$ of $\xi$  by the
formula
$$a(\xi) = \ds{\frac{\|\xi\|}{n} \sum_{k=2}^m \{\tr(A_1 A_k)\}\xi_k},$$
where $\{\xi_1 = \ds{\frac{\xi}{\|\xi\|}},\xi_2, \dots,  \xi_m \}$ is an
orthonormal base of the normal space of $M$, and $A_i = A_{\xi_i},
\,\, i = 1,\dots, m$ is the shape operator with respect to
$\xi_i$. The allied vector field $a(H)$ of the mean
curvature vector field $H$ is called the \emph{allied mean
curvature vector field} of $M$ in $\widetilde{M}$. B.-Y. Chen
defined  the $\mathcal{A}$-submanifolds to be those submanifolds
of $\widetilde{M}$ for which
 $a(H)$ vanishes identically \cite{Chen1}.
In \cite{GVV1}, \cite{GVV2} the $\mathcal{A}$-submanifolds are
called \emph{Chen submanifolds}. It is easy to see that minimal
submanifolds, pseudo-umbilical submanifolds and hypersurfaces are
Chen submanifolds. These Chen submanifolds are said to be trivial
Chen-submanifolds. In \cite{GM5} we showed that if  $M^2$ is a spacelike surface
in $\R^4_1$ with spacelike or timelike mean curvature vector field then the allied mean curvature vector field of $M^2$
is
$$a(H) = \ds{\frac{\sqrt{\varkappa^2-k}}{2} \,\lambda \, l}.$$
Hence, if $M^2$ is free of minimal points, then $a(H) = 0$ if and
only if $\lambda = 0$. This gives the geometric meaning of the
invariant $\lambda$:  $M^2$ is a non-trivial  Chen
surface if and only if the invariant $\lambda$ is zero.

\section{Basic classes of general rotational surfaces}

General rotational surfaces  in the Euclidean 4-space $\R^4$ were
introduced by  C. Moore \cite{M} as follows. Let $c: x(u) = \left(
x^1(u), x^2(u),  x^3(u), x^4(u)\right)$; $ u \in J \subset \R$ be
a smooth curve in $\R^4$, and $\alpha$, $\beta$ be constants. A
general rotation of the meridian curve $c$ in $\R^4$ is defined by
$$X(u,v)= \left( X^1(u,v), X^2(u,v),  X^3(u,v), X^4(u,v)\right),$$
where
$$\begin{array}{ll}
\vspace{2mm} X^1(u,v) = x^1(u)\cos\alpha v - x^2(u)\sin\alpha v; &
\qquad
X^3(u,v) = x^3(u)\cos\beta v - x^4(u)\sin\beta v; \\
\vspace{2mm} X^2(u,v) = x^1(u)\sin\alpha v + x^2(u)\cos\alpha v;&
\qquad X^4(u,v) = x^3(u)\sin\beta v + x^4(u)\cos\beta v.
\end{array}$$
In the case $\beta = 0$, $x^2(u) = 0$  the plane $Oe_3e_4$ is fixed and one gets
the classical rotation about a fixed two-dimensional axis.

In \cite{GM2-a} we considered  a special case of such surfaces,
given by
\begin{equation} \label{E:Eq-5}
\mathcal{M}: z(u,v) = \left(
f(u) \cos\alpha v, f(u) \sin \alpha v, g(u) \cos \beta v, g(u)
\sin \beta v \right),
\end{equation}
where $u \in J \subset \R, \,\,  v \in [0;
2\pi)$, $f(u)$ and $g(u)$ are smooth functions, satisfying
$\alpha^2 f^2(u)+ \beta^2 g^2(u)
> 0 , \,\, f'\,^2(u)+ g'\,^2(u) > 0$, and $\alpha, \beta$ are positive
constants. In the case  $\alpha \neq \beta$ each parametric curve
$u = const$ is a curve in $\R^4$ with constant Frenet curvatures,
and in the case  $\alpha = \beta$ each parametric curve $u =
const$ is a circle. The parametric curves  $v = const$  are plane curves with Frenet curvature $\ds{\frac{|g' f'' - f' g''|}{(\sqrt{f'\,^2 + g'\,^2})^3}}$.
These curves are the meridians of $\mathcal{M}$.

The surfaces defined by \eqref{E:Eq-5} are general rotational surfaces in
the sense of C. Moore with plane meridian curves. In \cite{GM2-a} we found the
invariants of these surfaces and completely  classified the minimal
super-conformal general rotational surfaces in $\R^4$. The classification of the general rotational surfaces in $\R^4$ consisting of parabolic points
is given in \cite{GM2-A}.

\vskip 2mm
Similarly to the general rotations in $\R^4$ one can consider
general rotational surfaces  in the Minkowski 4-space $\R^4_1$ as follows.
 Let $c: x(u) = \left( x^1(u), x^2(u),  x^3(u), x^4(u)\right)$; $ u \in J \subset \R$ be
a smooth spacelike or timelike curve in $\R^4_1$, and $\alpha$, $\beta$ be constants.
We consider the surface defined by
$$X(u,v)= \left( X^1(u,v), X^2(u,v),  X^3(u,v), X^4(u,v)\right),$$
where
$$\begin{array}{ll}
\vspace{2mm} X^1(u,v) = x^1(u)\cos\alpha v - x^2(u)\sin\alpha v; &
\qquad
X^3(u,v) = x^3(u)\cosh\beta v + x^4(u)\sinh\beta v; \\
\vspace{2mm} X^2(u,v) = x^1(u)\sin\alpha v + x^2(u)\cos\alpha v;&
\qquad X^4(u,v) = x^3(u)\sinh\beta v + x^4(u)\cosh\beta v.
\end{array}$$

In the case $\beta = 0$, $x^2(u) = 0$ (or  $x^1(u) = 0$)   one gets
the standard rotational surface of elliptic type in $\R^4_1$.
A local classification of spacelike rotational surfaces of elliptic type, whose
mean curvature vector field is either vanishing or lightlike, was obtained in \cite{Haesen-Ort-2}.

In the case $\alpha = 0$, $x^3(u) = 0$  we get
the standard hyperbolic rotational surface of first type, and in the case
$\alpha = 0$, $x^4(u) = 0$ we get the standard hyperbolic rotational surface of second type.
A local classification of spacelike rotational surfaces of hyperbolic type
with either vanishing or lightlike mean curvature vector field is given in \cite{Haesen-Ort-1}.
In \cite{Liu-Liu-1} the timelike and spacelike hyperbolic rotational surfaces
with non-zero constant mean curvature in the three-dimensional de Sitter space
$\mathbb{S}^3_1$ were classified.
Similarly, a classification of the spacelike and timelike Weingarten rotational surfaces in $\mathbb{S}^3_1$
is given in \cite{Liu-Liu-2}.
In \cite{GM3} we described the class of Chen spacelike rotational surfaces of hyperbolic or elliptic type in $\R^4_1$.

In the case $\alpha > 0$ and $\beta > 0$ the surfaces defined above are analogous to the general rotational surfaces of C. Moore in $\R^4$.

\vskip 2mm
In \cite{GM5} we considered  spacelike general rotational surfaces with plane meridians in the Minkowski space $\R^4_1$ and found their invariant functions.
Here we shall describe and classify some basic geometric classes of such surfaces.

 Let  $\mathcal{M}_1$ be the surface parameterized by
\begin{equation} \label{E:Eq-6}
\mathcal{M}_1: z(u,v) = \left( f(u) \cos \alpha v, f(u) \sin \alpha v, g(u) \cosh \beta v, g(u) \sinh \beta v \right),
\end{equation}
where $u \in J \subset \R$, $v \in [0; 2\pi)$, $f(u)$ and $g(u)$ are smooth functions, satisfying $\alpha^2 f^2(u)- \beta^2 g^2(u) > 0$, $f'\,^2(u)+ g'\,^2(u) > 0$,
and $\alpha, \beta$ are positive constants.

The coefficients of the first fundamental form of $\mathcal{M}_1$ are
$E = f'\,^2(u)+ g'\,^2(u)$;  $F = 0$; $G = \alpha^2 f^2(u)- \beta^2 g^2(u)$.
$\mathcal{M}_1$ is a spacelike surface whose mean curvature
vector at any point is a non-zero spacelike vector (see \cite{GM5}). Moreover, $\mathcal{M}_1$ is
 parameterized by principal parameters $(u,v)$.

The invariants $k$, $\varkappa$, and $K$ of
$\mathcal{M}_1$ are expressed by the functions $f(u)$, $g(u)$ and their derivatives as follows:
\begin{equation} \label{E:Eq-6-k}
k = \ds{\frac{4 \alpha^2 \beta^2 (g f' - f g')^2 (g' f'' - f' g'') (\alpha^2 f g' + \beta^2 g f')}{(f'\,^2 + g'\,^2)^3 (\alpha^2 f^2 - \beta^2 g^2)^3 }};
\end{equation}

\begin{equation} \label{E:Eq-6-kapa}
\varkappa =  \ds{\frac{\alpha \beta (g f' - f g') [(\alpha^2 f^2 - \beta^2  g^2)(g' f'' - f' g'') + (f'\,^2 +
g'\,^2) ( \alpha^2 f g' + \beta^2 g f') ]}{(f'\,^2 + g'\,^2)^2 (\alpha^2 f^2 - \beta^2 g^2)^2} };
\end{equation}

\begin{equation} \label{E:Eq-6-K}
K =  \ds{\frac{- (\alpha^2 f^2 - \beta^2  g^2)(\alpha^2 f g' + \beta^2 g f')(g' f'' - f' g'') + \alpha^2 \beta^2 (f'\,^2 + g'\,^2) (g f' - f g')^2}{(f'\,^2 + g'\,^2)^2 (\alpha^2 f^2 - \beta^2 g^2)^2 } \,}.
\end{equation}

The geometric invariant functions $\gamma_1, \, \gamma_2, \, \nu_1,\, \nu_2, \, \lambda, \, \mu,
\, \beta_1, \beta_2$ of $\mathcal{M}_1$ are:
\begin{equation} \label{E:Eq-7}
\begin{array}{ll}
\vspace{2mm}
\gamma_1 = 0; & \qquad \gamma_2 = \ds{- \frac{\alpha^2 f f' - \beta^2 g g'}{\sqrt{f'\,^2 + g'\,^2}(\alpha^2 f^2 - \beta^2 g^2)}};\\
\vspace{2mm} \nu_1 = \ds{\frac{g' f'' - f' g''}{(f'\,^2 +
g'\,^2)^{\frac{3}{2}}}}; & \qquad
\nu_2 = \ds{- \frac{\alpha^2 f g' + \beta^2 g f'}{\sqrt{f'\,^2 + g'\,^2}(\alpha^2 f^2 - \beta^2 g^2)}};\\
\vspace{2mm}
\lambda = 0; & \qquad  \mu = \ds{\frac{\alpha \beta (g f' - f g')}{\sqrt{f'\,^2 + g'\,^2}(\alpha^2 f^2 - \beta^2 g^2)}};\\
\vspace{2mm}
 \beta_1 = 0; & \qquad  \beta_2 = \ds{\frac{\alpha \beta (f f' + g g')}{\sqrt{f'\,^2 + g'\,^2}(\alpha^2 f^2 - \beta^2 g^2) }}.
\end{array}
\end{equation}

\vskip 3mm

In a similar way we can  consider the
surface $\mathcal{M}_2$ parameterized by
\begin{equation} \label{E:Eq-8}
\mathcal{M}_2: z(u,v) = \left( f(u) \cos \alpha v, f(u) \sin \alpha v, g(u) \sinh \beta v, g(u) \cosh \beta v \right),
\end{equation}
where $u \in J$, $v \in [0; 2\pi)$, $f(u)$ and $g(u)$ are smooth functions, satisfying
$f'\,^2(u)- g'\,^2(u) > 0$, $\alpha^2 f^2(u)+ \beta^2 g^2(u) > 0$,
and $\alpha, \beta$ are positive constants.

The coefficients of the first fundamental form
of $\mathcal{M}_2$ are
$E = f'\,^2(u)- g'\,^2(u)$; $F = 0$; $G =\alpha^2 f^2(u)+ \beta^2 g^2(u)$.
$\mathcal{M}_2$ is a spacelike surface with
timelike mean curvature vector field. The parameters $(u,v)$ of $\mathcal{M}_2$
are principal.

The invariants
 $k$, $\varkappa$, and $K$ of $\mathcal{M}_2$ are expressed similarly to the invariants of  $\mathcal{M}_1$  \cite{GM5}:
\begin{equation} \label{E:Eq-8-k}
k = \ds{\frac{4 \alpha^2 \beta^2 (g f' - f g')^2 (g' f'' - f' g'') (\alpha^2 f g' + \beta^2 g f')}{(f'\,^2 - g'\,^2)^3 (\alpha^2 f^2 + \beta^2 g^2)^3 }};
\end{equation}

\begin{equation} \label{E:Eq-8-kapa}
\varkappa =  \ds{\frac{\alpha \beta (g f' - f g')[(\alpha^2 f^2 + \beta^2  g^2)(g' f'' - f' g'') + (f'\,^2 -
g'\,^2) ( \alpha^2 f g' + \beta^2 g f') ]}{(f'\,^2 - g'\,^2)^2 (\alpha^2 f^2 + \beta^2 g^2)^2} };
\end{equation}

\begin{equation} \label{E:Eq-8-K}
K =  \ds{\frac{(\alpha^2 f^2 + \beta^2  g^2)(\alpha^2 f g' + \beta^2 g f')(g' f'' - f' g'') - \alpha^2 \beta^2 (f'\,^2 - g'\,^2) (g f' - f g')^2}{(f'\,^2 - g'\,^2)^2 (\alpha^2 f^2 + \beta^2 g^2)^2 } \,}.
\end{equation}

The geometric invariant functions  of
$\mathcal{M}_2$ are given below:
\begin{equation} \label{E:Eq-9}
\begin{array}{ll}
\vspace{2mm}
\gamma_1 = 0; & \qquad \gamma_2 = \ds{- \frac{\alpha^2 f f' + \beta^2 g g'}{\sqrt{f'\,^2 - g'\,^2}(\alpha^2 f^2 + \beta^2 g^2)}};\\
\vspace{2mm} \nu_1 = \ds{\frac{g' f'' - f' g''}{(f'\,^2 -
g'\,^2)^{\frac{3}{2}}}}; & \qquad
\nu_2 = \ds{- \frac{\alpha^2 f g' + \beta^2 g f'}{\sqrt{f'\,^2 - g'\,^2}(\alpha^2 f^2 + \beta^2 g^2)}};\\
\vspace{2mm}
\lambda = 0; & \qquad  \mu = \ds{\frac{\alpha \beta (f g' - g f')}{\sqrt{f'\,^2 - g'\,^2}(\alpha^2 f^2 + \beta^2 g^2)}};\\
\vspace{2mm}
 \beta_1 = 0; & \qquad  \beta_2 = \ds{\frac{\alpha \beta (g g' - f f')}{\sqrt{f'\,^2 - g'\,^2}(\alpha^2 f^2 + \beta^2 g^2) }}.
\end{array}
\end{equation}

We shall call the general rotational surface $\mathcal{M}_1$, defined by \eqref{E:Eq-6}, a \emph{general rotational surface of first type}, and
 the general rotational surface $\mathcal{M}_2$, defined by \eqref{E:Eq-8}, we shall call a \emph{general rotational surface of second type}.

Note that the invariant $\lambda$ of the general rotational surfaces of first or second type is zero. Hence, the following statement holds.

\begin{thm}
The  general rotational surfaces of first or second type,  free of minimal points, are non-trivial  Chen surfaces.
\end{thm}

In the following subsections we shall describe the classes of flat general rotational surfaces, general rotational surfaces with flat normal connection,
general rotational surfaces consisting of parabolic points, and minimal general rotational surfaces.

\subsection{Flat general rotational surfaces}

Let  $\mathcal{M}_1$ and $\mathcal{M}_2$  be  general rotational surfaces of first and second type, respectively.
Recall that a surface is called \emph{flat} if the Gauss curvature $K$ is zero.
Using equalities \eqref{E:Eq-6-K} and \eqref{E:Eq-8-K} we obtain

\begin{thm}\label{T:Th-flat}
(i) The general rotational surface of first type is flat if and only if
\begin{equation} \label{E:Eq-10}
\alpha^2 \beta^2 (f'\,^2 + g'\,^2) (g f' - f g')^2 = (\alpha^2 f^2 - \beta^2  g^2)(\alpha^2 f g' + \beta^2 g f')(g' f'' - f' g'').
\end{equation}

(ii) The general rotational surface of second type is flat if and only if
\begin{equation} \label{E:Eq-10-a}
\alpha^2 \beta^2 (f'\,^2 - g'\,^2) (g f' - f g')^2 = (\alpha^2 f^2 + \beta^2  g^2)(\alpha^2 f g' + \beta^2 g f')(g' f'' - f' g'').
\end{equation}
\end{thm}

Assume that the meridian curve is parameterized by $f = f(u); \,\, g = u$. Then  equation \eqref{E:Eq-10} takes the form
$$\frac{f''}{1 + f'\,^2} = \frac{\alpha^2 \beta^2 (u f' - f)^2 }{(\alpha^2 f^2 - \beta^2  u^2)(\alpha^2 f  + \beta^2 u f')},$$
which is equivalent to
\begin{equation} \label{E:Eq-11}
\left(\arctan f' \right)' = \frac{\alpha^2 \beta^2 (u f' - f)^2 }{(\alpha^2 f^2 - \beta^2  u^2)(\alpha^2 f  + \beta^2 u f')}.
\end{equation}

Similarly, equation \eqref{E:Eq-10-a} takes the form
$$\frac{f''}{1 - f'\,^2} = \frac{- \alpha^2 \beta^2 (u f' - f)^2 }{(\alpha^2 f^2 + \beta^2  u^2)(\alpha^2 f  + \beta^2 u f')},$$
which is equivalent to
\begin{equation} \label{E:Eq-11-a}
\left(\ln \left| \frac{1+f'}{1-f'} \right| \right)' = \frac{- 2\alpha^2 \beta^2 (u f' - f)^2 }{(\alpha^2 f^2 + \beta^2  u^2)(\alpha^2 f  + \beta^2 u f')}.
\end{equation}

Equations \eqref{E:Eq-11} and \eqref{E:Eq-11-a} describe analytically the class of flat general rotational surfaces of first and second type.

\subsection{General rotational surfaces with flat normal connection}
A surface is said to have  \emph{flat normal connection} if the curvature of the normal connection is zero.
The curvature of the normal connection of the general rotational surface  $\mathcal{M}_1$ (resp. $\mathcal{M}_2$) is
given by formula  \eqref{E:Eq-6-kapa} (resp. \eqref{E:Eq-8-kapa}). Using these formulas we obtain the next theorem.

\begin{thm}\label{T:Th-flat normal}
(i) The general rotational surface of first type has flat normal connection if and only if
\begin{equation} \label{E:Eq-12}
\frac{g' f'' - f' g''}{ f'\,^2 + g'\,^2} = - \frac{\alpha^2 f g' + \beta^2 g f'}{\alpha^2 f^2 - \beta^2  g^2}.
\end{equation}

(ii) The general rotational surface of second type has flat normal connection if and only if
\begin{equation} \label{E:Eq-12-a}
\frac{g' f'' - f' g''}{ f'\,^2 - g'\,^2} = - \frac{\alpha^2 f g' + \beta^2 g f'}{\alpha^2 f^2 + \beta^2  g^2}.
\end{equation}
\end{thm}

If we assume that  the meridian curve is parameterized by $f = f(u); \,\, g = u$, then  equation \eqref{E:Eq-12} takes the form
\begin{equation} \label{E:Eq-13}
\left(\arctan f' \right)' = - \frac{\alpha^2 f + \beta^2 u f'}{\alpha^2 f^2 - \beta^2  u^2}.
\end{equation}

Similarly, equation \eqref{E:Eq-12-a} takes the form
\begin{equation} \label{E:Eq-13-a}
\left(\ln \left| \frac{1+f'}{1-f'} \right| \right)' = \frac{ 2(\alpha^2 f + \beta^2 u f') }{\alpha^2 f^2 + \beta^2  u^2}.
\end{equation}

The class of  general rotational surfaces  with flat normal connection is described analytically  by equations \eqref{E:Eq-13} and \eqref{E:Eq-13-a}.

\vskip 3mm \noindent
\textbf{Example 1.}
Let $f(u) = a \cos u,\,\, g(u) = a \sin u$, $a = const$ ($a \neq 0$). A direct computation shows that equation \eqref{E:Eq-12} is fulfilled.
Hence, the surface parameterized by
$$z(u,v) = \left(a \cos u \cos \alpha v, a \cos u \sin \alpha v, a \sin u \cosh \beta v, a \sin u \sinh \beta v \right)$$
is a spacelike general rotational surface of first type with flat normal connection.

In the special case when $a = 1$,
$\alpha=\beta=1$ we obtain a spacelike surface lying on the De Sitter
space $S^3_1 = \{x \in \R^4_1; \langle x,x \rangle = 1\}$.

\vskip 3mm \noindent
\textbf{Example 2.}
If we choose   $f(u) = a \sinh u,\,\, g(u) = a \cosh u$, $a = const$ ($a \neq 0$), by a direct computation we obtain that equation \eqref{E:Eq-12-a} is fulfilled.
Hence, the surface parameterized by
$$z(u,v) = \left( a \sinh u \cos \alpha v, a \sinh u \sin \alpha v, a \cosh u \sinh \beta v, a \cosh u \cosh \beta v \right)$$
is a spacelike general rotational surface of second type with flat normal connection.

In the special case  $a = 1$,
$\alpha=\beta=1$ we obtain a spacelike surface lying on the unit
hyperbolic sphere $H^3_1 = \{x \in \R^4_1; \langle x,x \rangle = - 1\}$.

\subsection{General rotational surfaces consisting of parabolic points}

Recall that surfaces consisting of parabolic points are characterized by the condition $k = 0$.
The next theorem classifies the general rotational surfaces of first and second type with $k = 0$.

\begin{thm}\label{T:Th-parabolic points}
A general rotational surface of first or second type consists of parabolic points if and only if it is one of the following:

(i) a developable ruled surface in $\R^4_1$;

(ii) a non-developable ruled surface in $\R^4_1$;

(iii) a non-ruled surface in $\R^4_1$ whose meridian curve is given by $f(u) = u; \,\, g(u) = \ds{ c\,u^{-\frac{\beta^2}{\alpha^2}}}$, $c = const \neq 0$.
\end{thm}

\vskip 2mm \noindent
\emph{Proof:}
Consider a general rotational surface of first or second type.
 Equality \eqref{E:Eq-6-k} (or \eqref{E:Eq-8-k}) implies that $k = 0$ if and only if
\begin{equation} \label{E:Eq-14}
(g f' - f g') (g' f'' - f' g'') (\alpha^2 f g' + \beta^2 g f') = 0.
\end{equation}
Without loss of generality we assume that the
meridian curve is given by $f = u; \,\, g = g(u)$.
Then equality  \eqref{E:Eq-14} takes the form
\begin{equation} \notag
(g - u g')  g'' (\alpha^2 u g' + \beta^2 g) = 0,
\end{equation}
which implies that the invariant $k$ is zero in the following three cases:

\vskip 1mm 1. $g(u) = a\,u,\; a = const \neq 0$. In this case $k =
\varkappa = K = 0$, and by a result in \cite{GM5} the corresponding general rotational surface (of first or second type) is a developable ruled surface in $\R^4_1$.

\vskip 1mm 2. $g(u) = a\,u + b ,\; a= const \neq 0, b = const \neq
0$.  Hence, the meridians are straight lines. It can easily be seen that in this case $\varkappa \neq 0$.
Consequently, the corresponding general rotational surface is a non-developable ruled surface in $\R^4_1$.

\vskip 1mm 3. $g(u) = \ds{ c\,u^{-\frac{\beta^2}{\alpha^2}}}, \; c = const \neq 0$.
In  this case the meridians are not straight lines. The invariants $\varkappa$ and $K$ are non-zero, and hence, the
corresponding general rotational surface is a non-ruled surface in $\R^4_1$.

\qed

\subsection{Minimal general rotational surfaces}
In this subsection we shall find all minimal  general rotational surfaces of first and second type.
Recall that a surface is minimal if and only if $\nu_1 + \nu_2 = 0$.
Hence, using \eqref{E:Eq-7} we get that
the general rotational surface of first type is minimal if and only if the functions $f(u)$ and $g(u)$ satisfy the equality
\begin{equation} \label{E:Eq-15}
\frac{g' f'' - f' g''}{ f'\,^2 + g'\,^2} =  \frac{\alpha^2 f g' + \beta^2 g f'}{\alpha^2 f^2 - \beta^2  g^2}.
\end{equation}
Similarly, from \eqref{E:Eq-9} it follows that the general rotational surface of second type is minimal if and only if
\begin{equation} \label{E:Eq-15-a}
\frac{g' f'' - f' g''}{ f'\,^2 - g'\,^2} =  \frac{\alpha^2 f g' + \beta^2 g f'}{\alpha^2 f^2 + \beta^2  g^2}.
\end{equation}

We shall find the solutions of equalities \eqref{E:Eq-15} and \eqref{E:Eq-15-a}. In such a way
we shall describe the class of minimal  general rotational surfaces of first and second type.

\begin{thm}\label{T:Th-minimal}
(i) The general rotational surface of first type is minimal if and only if
the meridian curve is given by the formula
$$g = \frac{\sqrt{A}}{\beta} \sin \left( \varepsilon \frac{\beta}{\alpha} \ln \left| \alpha f + \sqrt{\alpha^2 f^2 - A} \right| + C \right), \qquad C = const, \,\, A = const > 0.$$

(ii) The general rotational surface of second type is minimal if and only if
the meridian curve is given by the formula
$$g = \frac{\sqrt{A}}{\beta} \sin \left( \varepsilon \frac{\beta}{\alpha} \ln \left| \alpha f + \sqrt{\alpha^2 f^2 + A} \right| + C \right), \qquad C = const, \,\, A = const >0.$$
\end{thm}

\vskip 2mm \noindent
\emph{Proof:}
First, we shall simplify equalities \eqref{E:Eq-15} and \eqref{E:Eq-15-a}.
Using that the connection $\nabla'$ of $\R^4_1$  is flat, from $R'(x,y,x) = 0$ and  $R'(x,y,y) = 0$ we obtain that
 the invariants of each spacelike surface in $\R^4_1$ satisfy the equalities (see \cite{GM5}):
\begin{equation} \label{E:Eq-16}
\begin{array}{l}
\vspace{2mm}
2\mu\, \gamma_2 + \nu_1\,\beta_2 - \lambda\,\beta_1 = x(\mu);\\
\vspace{2mm}
2\lambda\, \gamma_1 - \mu\,\beta_2 + (\nu_1 - \nu_2)\,\gamma_2 = - x(\nu_2) + y(\lambda).
\end{array}
\end{equation}
In the case the surface is a minimal general rotational surface of first or second type we have $\lambda = 0$, $\gamma_1 = 0$, $\beta_1 = 0$ and
$\nu_2 = - \nu_1$. Hence, from \eqref{E:Eq-16} it follows that
\begin{equation} \notag
\begin{array}{l}
\vspace{2mm}
2\mu\, \gamma_2 + \nu_1\,\beta_2 = x(\mu);\\
\vspace{2mm}
2\nu_1\, \gamma_2 - \mu\,\beta_2 =  x(\nu_1),
\end{array}
\end{equation}
which imply
$ \gamma_2 = \ds{\frac{1}{4}} \,x\left( \ln(\mu^2 + \nu_1^2)\right)$.
On the other hand,  $\gamma_2 = - x\left( \ln \sqrt{G}\right)$.
Hence, we get
$\ds{\frac{1}{4}} \,x\left( \ln(\mu^2 + \nu_1^2)\right) + x\left( \ln \sqrt{G}\right))=  0$,
which implies
 $$x\left( G^2(\mu^2 + \nu_1^2)\right) = 0.$$
Now, using that
$\mu$, $\nu_1$, and $G$ are functions depending only on the parameter $u$, we
obtain
\begin{equation}  \label{E:Eq-17}
 G^2 (\mu^2 + \nu_1^2) = c^2,
\end{equation}
where $c$ is a constant.

\vskip 2mm
Now, let $\mathcal{M}_1$  be a minimal general rotational surface of first type. Then $G = \alpha^2 f^2 - \beta^2 g^2$, and using  \eqref{E:Eq-7} and \eqref{E:Eq-17} we obtain

\begin{equation} \notag
\frac{\alpha^2 \beta^2 (g f' - f g')^2 + (\alpha^2 f g' + \beta^2 g f')^2}{f'\,^2 + g'\,^2} = c^2,
\end{equation}
which is equivalent to
\begin{equation} \label{E:Eq-18}
\frac{\alpha^2 f^2 g'\,^2 + \beta^2 g^2 f'\,^2}{f'\,^2 + g'\,^2} = \frac{c^2}{\alpha^2 + \beta^2}.
\end{equation}
Equality \eqref{E:Eq-18} can  be obtained also from \eqref{E:Eq-15} by  a direct but very long computation.

Without loss of generality we  assume that $f'\,^2 + g'\,^2 = 1$. Then from \eqref{E:Eq-18} we get
\begin{equation} \label{E:Eq-19}
\alpha^2 f^2 g'\,^2 + \beta^2 g^2 f'\,^2 = \frac{c^2}{\alpha^2 + \beta^2}.
\end{equation}
Denote $A = \ds{\frac{c^2}{\alpha^2 + \beta^2}}$.
Now, using that $g'\,^2 = 1 - f'\,^2$, from \eqref{E:Eq-19} it follows that
\begin{equation} \label{E:Eq-20}
f'\,^2 = \frac{\alpha^2 f^2 - A}{\alpha^2f^2 - \beta^2 g^2}; \qquad g'\,^2 = \frac{A - \beta^2 g^2}{\alpha^2f^2 - \beta^2 g^2}.
\end{equation}
Note that the constant $A$ satisfies $\beta^2 g^2 < A < \alpha^2 f^2$, since $\alpha^2f^2 - \beta^2 g^2 > 0$.

Equalities \eqref{E:Eq-20} imply $(A - \beta^2 g^2) f'\,^2 = (\alpha^2 f^2 - A) g'\,^2$, i.e.
$$\frac{f'}{\sqrt{\alpha^2 f^2 - A}} = \varepsilon \frac{g'}{\sqrt{A - \beta^2 g^2}}, \qquad \varepsilon =\pm 1.$$
Integrating the last equality we obtain
$$\int \frac{df}{\sqrt{\alpha^2 f^2 - A}} = \varepsilon \int \frac{dg}{\sqrt{A - \beta^2 g^2}}.$$
Calculating the integrals we get
$$\arcsin \frac{\beta g}{\sqrt{A}} =\varepsilon \frac{\beta}{\alpha} \ln \left| \alpha f + \sqrt{\alpha^2 f^2 - A} \right| + C, \qquad C = const.$$

Consequently, in the case of a minimal general rotational surface of first type the meridian curve is given by the following formula
\begin{equation} \label{E:Eq-21}
g = \frac{\sqrt{A}}{\beta} \sin \left( \varepsilon \frac{\beta}{\alpha} \ln \left| \alpha f + \sqrt{\alpha^2 f^2 - A} \right| + C \right).
\end{equation}

Conversely, if the meridian curve is defined by  formula \eqref{E:Eq-21}, by a straightforward computation we obtain that equality \eqref{E:Eq-15} is fulfilled
and hence, the general rotational surface of first type is minimal.

\vskip 2mm
In the case of  a minimal general rotational surface of second type we have $G = \alpha^2 f^2 + \beta^2 g^2$. Then  \eqref{E:Eq-9} and \eqref{E:Eq-17} imply
\begin{equation} \label{E:Eq-18-a}
\frac{\alpha^2 f^2 g'\,^2 + \beta^2 g^2 f'\,^2}{f'\,^2 - g'\,^2} = \frac{c^2}{\alpha^2 + \beta^2}.
\end{equation}

Without loss of generality we  assume that $f'\,^2 - g'\,^2 = 1$. Then from \eqref{E:Eq-18-a} we get
\begin{equation} \label{E:Eq-19-a}
\alpha^2 f^2 g'\,^2 + \beta^2 g^2 f'\,^2 = A,
\end{equation}
where $A = \ds{\frac{c^2}{\alpha^2 + \beta^2}}$.
Using that $g'\,^2 = f'\,^2 - 1$, from \eqref{E:Eq-19-a} we obtain
\begin{equation} \label{E:Eq-20-a}
f'\,^2 = \frac{A + \alpha^2 f^2}{\alpha^2f^2 + \beta^2 g^2}; \qquad g'\,^2 = \frac{A - \beta^2 g^2}{\alpha^2f^2 + \beta^2 g^2}.
\end{equation}
Note that in this case $A > \beta^2 g^2$, since $g'\,^2 > 0$.

Equalities \eqref{E:Eq-20-a} imply
$$\frac{f'}{\sqrt{A + \alpha^2 f^2}} = \varepsilon \frac{g'}{\sqrt{A - \beta^2 g^2}}, \qquad \varepsilon =\pm 1.$$
After integration we obtain
$$\arcsin \frac{\beta g}{\sqrt{A}} =\varepsilon \frac{\beta}{\alpha} \ln \left| \alpha f + \sqrt{A + \alpha^2 f^2} \right| + C, \qquad C = const.$$

Consequently, in the case of a minimal general rotational surface of second type the meridian curve is given by
\begin{equation} \label{E:Eq-21-a}
g = \frac{\sqrt{A}}{\beta} \sin \left( \varepsilon \frac{\beta}{\alpha} \ln \left| \alpha f + \sqrt{A + \alpha^2 f^2} \right| + C \right).
\end{equation}

A direct computation shows that if the meridian curve is given by formula \eqref{E:Eq-21-a}, then equality \eqref{E:Eq-15-a} is satisfied.
Hence, the  general rotational surface of second type is minimal.

\qed

\vskip 2mm
Finally, it should be mentioned that the classes of minimal general rotational surfaces
and  general rotational surfaces consisting of parabolic points are found explicitly. The  classes of flat  general rotational surfaces and
 general rotational surfaces with flat normal connection are described analytically by ordinary differential equations.
 An open question is to find them explicitly.

\end{document}